\documentclass[a4paper,12pt]{amsart}

\usepackage{amssymb}
\usepackage{amsmath}

\newtheorem{thm}[equation]{Theorem}
\newtheorem{lem}[equation]{Lemma}
\newtheorem{defn}{Definition}[section]

\usepackage{ifthen}
\usepackage{setspace}
\newenvironment{pf}[1][]{%
 \vskip 3mm
 \noindent
 \ifthenelse{\equal{#1}{}}%
  {{\slshape Proof. }}%
  {{\slshape #1.} }%
 }%
{\qed\bigskip}

\begin{document}
\bibliographystyle{amsplain}

\newcommand{\A}{{\mathcal A}}
\newcommand{\D}{{\mathbb D}}
\newcommand{\es}{{\mathcal S}}
\newcommand{\CC}{{\mathcal C}}
\newcommand{\K}{{\mathcal K}}
\newcommand{\IC}{{\mathbb C}}
\newcommand{\blem}{\begin{lem}}
\newcommand{\elem}{\end{lem}}
\newcommand{\bthm}{\begin{thm}}
\newcommand{\ethm}{\end{thm}}
\newcommand{\bpf}{\begin{pf}}
\newcommand{\epf}{\end{pf}}
\newcommand{\bdefe}{\begin{defn}}
\newcommand{\edefe}{\end{defn}}

\title[Close-to-Convexity properties of Clausen's  Hypergeometric Function $_3F_2(a,b,c;d,e;z)$]{Close-to-Convexity properties of Clausen's  Hypergeometric Function $_3F_2(a,b,c;d,e;z)$}

\author{K.Chandrasekran} \address{K. Chandrasekran, Department of Mathematics, Jeppiaar SRR Engineering College, Padur, Chennai- 603 103, India.}\email{kchandru2014@gmail.com}

\author{D.J. Prabhakaran} \address{D.J. Prabhakaran, Department of Mathematics, Anna University, MIT Campus , Chennai- 600 005, India.}\email{asirprabha@yahoo.com}


\keywords{Hypergeometric Function, Starlike, Univalent, Convex, and Close-to-Convex Functions. \\ \hfill \textbf{Final Version after Correction as on 31.10.2020}}


\begin{abstract}
In this paper,  we obtain various conditions on the parameters $a,\, b,\, c\,, d$ and $e$ for which the hypergeometric functions  $z\, _3F_2(a,b,c;d,e;z)$ to be in the class of all close-to-convex function with respect to some well known convex functions.
\end{abstract}

\thanks{}

\maketitle

\markboth{K.CHANDRASEKRAN AND D.J. PRABHAKARAN }
{Close-to-convexity properties of $_3F_2(a,b,c;d,e;z)$}
\pagestyle{myheadings}

\section{Introduction and preliminaries}
Let $\A$ denote the family of functions $f$ that are analytic in the unit disk $\D =\{z:\, |z|<1\}$ and normalized by $f(0)=f'(0)-1=0$ and let $\es\subset \A $ be the class of all univalent functions. In  \cite{Duren-1983-book, A-W-Goodman-1983-book}, Various subclass of univalent functions are characterised by its geometrical properties and  the important subclasses are convex, starlike and close-to-convex functions.\\
\bdefe
A function $f$ in $\es$ is a convex function   in $\D$ if and only if
 \begin{eqnarray*}
  { \rm Re }\left( 1 + \frac{zf^{\prime\prime}(z)}{f^\prime(z)} \right) > 0
 \end{eqnarray*} Let us denote by ${\CC}$ the class of all convex univalent functions in $\mathbb{D}$.\\
\edefe
\bdefe
A function $f$ in $\es$ is a starlike function with respect to the origin in $\D$ if and only if
 \begin{eqnarray*}
{ \rm Re }  \left( \frac{zf^{\prime}(z)}{f(z)} \right) > 0
 \end{eqnarray*}
\edefe
Let us denote by $\es^{\ast}$, the class of all starlike univalent functions in $\D$.
\bdefe
A function $f$ in $\es$ is a close to convex in $\D$ if and only if
 \begin{eqnarray*}
{ \rm Re }  \left( \frac{f^{\prime}(z)}{g^{\prime}(z)} \right) > 0
 \end{eqnarray*}
for some $ g \in \CC $.
\edefe

The class of all close to convex univalent functions in $\D$ is denoted by $\K$. It is well known that the chain of inclusion relations $\CC \subset \es^{*} \subset \K \subset \es$.\\

The family of functions in $\A$ which are close-to-convex with respect to $-\log (1-z)$, and also, starlike in $\D$ is denoted by  $KS^{*}$.\\

\bdefe
Let $ \displaystyle f(z)= z+\sum_{n=2}^{\infty}\, a_n\,z^n $ and  $ \displaystyle g(z)= z+\sum_{n=2}^{\infty}\, b_n\,z^n $ be analytic in $\D$. Then the Hadamard product or convolution  of $f(z)$  and $g(z)$ is defined by
 \begin{eqnarray*}
\displaystyle f(z)*g(z)= z+\sum_{n=2}^{\infty} a_nb_n z^n,\, |z|<1.
 \end{eqnarray*}
\edefe
\bdefe
The Alexander transform, is defined as
 \begin{eqnarray}\label{alexeq1}
\Lambda_f(z) = \int_{0}^{z}\frac{f(t)}{t}\ dt, \ f \in \es , z \in \D
 \end{eqnarray}
\edefe
Using convolution technique the above transform is given by $ F(z) = f(z) \ast h(z) $
where $h(z) = - \log \left( 1-z \right)$, is not univalent in $\D$. It is know that the function $h(z)$ is well known and is a convex function of order $1/2$.\\

For any complex variable $a\neq 0$, the shifted factorial (or Pochhammer symbol) is defined as
$$(a)_0\,=\,1,\quad (a)_n\,=\,a(a+1)\cdots (a+n-1),\quad n\,=\,1,2,3,\cdots$$

In terms of  Euler gamma function,  the Pochhammer symbol can also be defined as $$(a)_n = \frac{\Gamma(n+a)}{\Gamma(a)},\quad n=0,1,2,\cdots$$
where $a$ is neither zero nor a negative integer.\\

\bdefe
The Hypergeometric Function $_3F_2(a,b,c;d,e;z)$ is defined by
 \begin{eqnarray}\label{inteq5}
_3F_2(a,b, c;d,e;z)=\sum_{n=0}^{\infty}\frac{(a)_n(b)_n(c)_n}{(d)_n(e)_n(1)_n}z^n;\, \, \, a,b,c,d,e\in \IC
 \end{eqnarray}
\edefe
where  $ a,b,c,d,e \in \IC$ and $d,e \neq 0, -1, -2, -3, \cdots$ which is analytic and convergent in  unit disc $\D$.\\

Using (\ref{inteq5}), the normalized hypergeometric function be written as $z\, _3F_2(a,b,c;d,e;z)$ and defined by
 \begin{eqnarray}\label{eq1}
f(z)=z\, _3F_2(a,b,c;d,e;z)=z+\sum_{n=2}^{\infty} A_n z^n,
 \end{eqnarray}
where
 \begin{eqnarray} \label{2f2eq1}
A_n = \frac{(a)_{n-1}(b)_{n-1}(c)_{n-1}}{(d)_{n-1}(e)_{n-1}(1)_{n-1}}, for\ n\geq 2.
 \end{eqnarray}
and $A_1=1$.

\blem \label{lemeq5} \cite{Fejer-1936-Trans-ams}
 If $A_n \geq 0, \{ nA_n\}$ and $\{nA_n-(n+1)A_{n+1}\}$ both are non-increasing, i.e., $\{nA_n\}$ is monotone of order 2, then  $f$ defined by (\ref{eq1})  is in $\es^{*}.$
\elem

\blem \label{lemeq2} \cite{Ozaki-1935-Tokyo} Suppose that
 \begin{eqnarray*} \label{lmeq1}
1\geq 2 A_2 \geq \cdots \geq nA_n \geq \cdots \geq 0
 \end{eqnarray*}
or
 \begin{eqnarray*} \label{lmeq2}
1\leqq 2 A_2 \leq \cdots \leq nA_n \leq \cdots \leq 2
 \end{eqnarray*}
Then $f(z)$ is defined by (\ref{eq1}) is close-to-convex with respect to $-\log(1-z)$.
\elem
\blem \label{lemeq3} \cite{Ozaki-1935-Tokyo}
Suppose that $f$ is an odd function ( i.e., the value of $A_{2n}$ in $(\ref{eq1})$ is zero for each $n\geq1$ ) such that
 \begin{eqnarray*} \label{lmeq3}
1\geq 3 A_3 \geq \cdots \geq (2n+1)A_{2n+1} \geq \cdots \geq 0
 \end{eqnarray*}
or
 \begin{eqnarray*} \label{lmeq4}
1\leq 3 A_3 \leq \cdots \leq (2n+1)A_{2n+1} \leq \cdots \leq 2
 \end{eqnarray*} %
Then $f\in \es$. In fact, $f(z)$ is close-to-convex with respect to the convex function $ \frac 1 2 \log((1+z)/(1-z))$.
\elem
In 1986, Ruscheweyh and Singh \cite{Ruscheweyh-and-Singh-1986} obtained the sufficient conditions on the parameters $a,\,b\,$ and $c$ for $z\,_2F_1(a,b;c;z)$ to be starlike of order $\beta < 1$. Further the year 1995, Ponnusamy and Vourinen \cite{Ponnusamy-Vuorinen-1995} have established the univalence and convexity properties of Gaussian hypergeometric functions. Subsequently Ponnusamy derived the conditions on the parameters $a,\,b$ and $c$ for univalence and starlike properties of Alexander transform and also the close-to-convexity properties of the Gaussian hypergeometric functions \cite{Ponnusamy-1996,Ponnusamy-1998}.\\

Inspired by the above results, in this paper, we find conditions on $ a, b, c, d$ and $e$ such that the function $ z\ _3F_2(a,b,c;d,e;z) $
to be close-to-convex with respect to the functions  $ -\log(1-z),\,  \frac{1}{2}\log\left(\frac{1+z}{1-z}\right)$ and is in the class $ KS^{\ast}$. We also find similar conditions such that the Alexander transform is in the class $KS^{\ast}$.

\section{Main Results and Proofs:}

\bthm\label{thm12f2}
If $a,b,c > 0$,  $de \geq 2\,a\,b\,c\, $  and
 \begin{eqnarray*}\label{thm12f2eq1}
d+e \geq Max\left\{ a+b+c,\, \frac{1}{2}(ab+bc+ac+2(a+b+c)-1-2abc),2[ab+bc+ac]-3abc \right\},
 \end{eqnarray*} then $z\ _3F_2(a,b,c;d,e;z)$ is close-to-convex with respect to $-log(1-z)$.
\ethm

\bpf Consider the function $f(z)$ is defined in $(\ref{eq1})$. Replace $n$ by $n+1$ in the equation (\ref{2f2eq1}) and Using the Pochhammer symbol, we have
 \begin{eqnarray}\label{2f2eq2}
A_{n+1} &=& \frac{(a)_{n}(b)_{n}(c)_{n}}{(d)_{n}(e)_{n}(1)_{n}} = \frac{(a+n-1)(b+n-1)(c+n-1)}{(d+n-1)(e+n-1)n} A_n
 \end{eqnarray}
and observe that $A_n >0$ for all $n\geq 1$. \\

To prove $f$ is close-to-convex with respect to $-\log(1-z)$. It is enough to show that $\{nA_n\}$ is non-increasing sequence.\\

From (\ref{2f2eq1}) and (\ref{2f2eq2}), we have the following after some manipulation
 \begin{eqnarray*}
nA_n - (n+1)A_{n+1} &=& nA_n -  \frac{(n+1)(a+n-1)(b+n-1)(c+n-1)}{(d+n-1)(e+n-1)n} A_n\\
&=& \frac{A_n U(n) }{(d+n-1)(e+n-1)n}
 \end{eqnarray*}
where
 \begin{eqnarray*}\label{2f2eq3}
U(n) &=&n^2(e+n-1)(d+n-1)-(n+1)(a+n-1)(b+n-1)(c+n-1)\nonumber\\
&=& {\left( e+d-c-b-a\right)} \,{n}^{3}+\left(d\,e-(e+d)-(a\,b+b\,c+a\,c)+c+b+a+1\right) \,{n}^{2}\nonumber\\ && -\left( a\,b\,c-c-b-a+2\right) \,n-\left( a-1\right) \,\left( b-1\right) \,\left( c-1\right)
 \end{eqnarray*}
For every $n\geq 1$, we have $n^3 \geq 3\,n^2-3\,n+1,$
 \begin{eqnarray*}
V(n) &\geq & {\left(d\,e+2\,(e+d)-(a\,b+b\,c+a\,c)-2\,(a+b+c)+1\right)} \,{n}^{2}\\
&& \qquad +\left(-3\,(e+d)-a\,b\,c+4\,(c+b+a)-2\right) \,n\\
&& \qquad \qquad +e+d-a\,b\,c+b\,c+a\,c-2\,c+a\,b-2\,b-2\,a+1
 \end{eqnarray*}
Using the fact that $n^2 \geq 2n-1$, for all $n\geq1.$
 \begin{eqnarray*}
W(n) &\geq & {\left(e+d+2\,d\,e-a\,b\,c-2\,(a\,b+b\,c+a\,c)\right)} \,n\\
        && \qquad  -\left( d\,e+e+d+a\,b\,c-2\,(a\,b+b\,c+\,a\,c)\right)
 \end{eqnarray*}
Put $n=1$ in above the equation, we find
 \begin{eqnarray*}
W(1) &=&{ d\,e-2\,a\,b\,c} \geq 0
 \end{eqnarray*}
Since $d\,e \geq 2\,a\,b\,c$, The above  equation implies that
$$ U(n) \geq  V(n) \geq  W(n) \geq W(1) \geq 0.$$
Hence $\{n A_n\}$ is a non-increasing sequence). Hence by Lemma  \ref{lemeq2}, the function \\$z\ _3F_2(a,b,c;d,e;z)$ is close-to-convex with respect to $-log(1-z)$.
\epf

\bthm\label{thm22f2}
If $a,b,c > 0$,  $de \geq 3abc $ and
 \begin{eqnarray}\label{thm22f2eq1}
d+e \geq Max\left\{a+b+c, \alpha(a,b,c),\, 3(a\,b+b\,c+a\,c)-7a\,b\,c\right\}
 \end{eqnarray}
where $$\alpha(a,b,c)=\frac{1}{3}\left((2(a\,b+b\,c+a\,c)+3(a+b+c)-6\,a\,b\,c-1\right)$$
then $z\ _3F_2(a,b,c;d,e;z^2)$ is close-to-convex with respect to $\frac{1}{2}\log((1+z)/(1-z))$.
\ethm
\bpf Consider the function defined as follows by replacing $z$ by $z^2$ in the equation (\ref{eq1})
 \begin{eqnarray*}
f(z)=z\ _3F_2(a,b,c;d,e;z^2)=z+\sum_{n=2}^{\infty} A_{2n-1} z^{2n-1},
 \end{eqnarray*}
where
 \begin{eqnarray}\label{2f2eq4}
A_{2n-1} = \frac{(a)_{n-1}(b)_{n-1}(c)_{n-1}}{(d)_{n-1}(e)_{n-1}(1)_{n-1}}, for\ n\geq 2.
 \end{eqnarray}
and $A_1=1$. Then,  we have the following from the equation (\ref{2f2eq4}) by replacing $n$ by $n+1$.
 \begin{eqnarray}\label{2f2eq6}
A_{2n+1} &=&  \frac{(a+n-1)(b+n-1)(c+n-1)}{(d+n-1)(e+n-1)n} A_{2n-1}\nonumber
 \end{eqnarray}
and therefore we obtain
 \begin{eqnarray*}
(2n-1)A_{2n-1} - (2n+1)A_{2n+1} &=&  \frac{A_{2n-1}  X(n)}{(e+n-1)(d+n-1)n}
 \end{eqnarray*}
where
 \begin{eqnarray*}
X(n) &=&2{\left( e+\,d-\,c-\,b-\,a\right)} \,{n}^{3}\\
        && \qquad +\left(2\,d\,e-3\,(e+d)-2\,(a\,b+b\,c+\,a\,c)+3\,(c+b+a)+1\right) \,{n}^{2}\\
        && \qquad \qquad \qquad -\left( d\,e-(e+d)+2\,a\,b\,c-(a\,b+b\,c+a\,c)+2\right) \,n \\
        && \qquad \qquad \qquad \qquad \qquad -\left( a-1\right) \,\left( b-1\right) \,\left( c-1\right)
 \end{eqnarray*}
Replace $n^3 \geq 3\,n^2-3\,n+1,$ for every $n\geq 1$, we have
 \begin{eqnarray*}
Y(n) &=& {\left(2\,d\,e+3\,(e+d)-2\,(a\,b+b\,c+\,a\,c)-3\,(c+b+a)+1\right)}\,n^2\\
        &&\qquad -\left( d\,e+5\,(e+d)+2\,a\,b\,c-(b\,c+a\,c+a\,b)-6\,(c+\,b+\,a)+2\right) \,n \\
                 &&\qquad \qquad\qquad + 2\,(e+d)-a\,b\,c+a\,b+b\,c+a\,c-3\,(c+b+a)+1
 \end{eqnarray*}
Using  $n^2 \geq 2n-1$, for all $n\geq1$, implies that
 \begin{eqnarray*}
Z(n) &=& {\left(3\,d\,e+e+d-2\,a\,b\,c-3(a\,b+\,b\,c+a\,c)\right)}  \,n\\
            && \qquad -\left( 2\,d\,e+e+d+a\,b\,c-3(a\,b+\,b\,c+a\,c)\right)
 \end{eqnarray*}
Replace $n$ by 1 in above, we get
 \begin{eqnarray*}
Z(1) &=& {d\,e-3\,a\,b\,c }\geq 0
 \end{eqnarray*}
Since $d\,e \geq 3\,a\,b\,c$, The above  equation implies that
$$ X(n) \geq  Y(n) \geq Z(n)\geq Z(1) \geq 0$$
Hence by the condition on $d+e$ in (\ref{thm22f2eq1}), we have $X(n)$ is non-negative, for all $n\geq 1$. Thus $\{(2n-1)A_{2n-1}\}$ is non-increasing sequence. The function $z\ _3F_2(a,b,c;d,e;z^2)$ is close-to-convex with respect to $\frac{1}{2}\log((1+z)/(1-z))$ by Lemma \ref{lemeq3}.
\epf

\bthm\label{thm32f2}
Let $a,\, b,\,and\, c > 0$, $$d+e \geq Max\left\{ T_1(a,b,c),\,  T_2(a,b,c),\,  T_3(a,b,c),\, T_4(a,b,c) \right\}$$
where
\begin{eqnarray*}
T_1(a,b,c) &=& \left( e+d-c-b-a\right) \,\left( e+d-c-b-a+1\right)\\
T_2(a,b,c) &=& \left( e+d-c-b-a+1\right) \,\left( 2\,d\,e+5\,(e+d)-2\,(a\,b+b\,c+a\,c)-5\,(c+b+a)+2\right)\\
T_3(a,b,c) &=& \left(\left( {d}^{2}+9\,d+9\right) \,{e}^{2}+\left( 9\,{d}^{2}+\left( \left( -2\,b-2\,a-8\right) \,c+\left( -2\,a-8\right) \,b-8\,a+29\right) \,d \right.\right.\\
&& \qquad \left.\left.+\left( \left( -2\,a-10\right) \,b-10\,a-16\right) \,c+\left( -10\,a-16\right) \,b-16\,a+15\right) \,e+9\,{d}^{2} \right.\\
&& \qquad \qquad \left.+\left( \left( \left( -2\,a-10\right) \,b-10\,a-16\right) \,c+\left( -10\,a-16\right) \,b-16\,a+15\right) \,d\right.\\
&& \qquad \qquad \qquad \left.+\left( {b}^{2}+\left( 4\,a+9\right) \,b+{a}^{2}+9\,a+7\right) \,{c}^{2}+\left( \left( 4\,a+9\right) \,{b}^{2}\right.\right.\\
&& \qquad \qquad \qquad \qquad \left.\left.+\left( 4\,{a}^{2}+24\,a+3\right) \,b+9\,{a}^{2}+3\,a-11\right) \,c+\left( {a}^{2}+9\,a+7\right) \,{b}^{2}\right.\\
&& \qquad \qquad \qquad \qquad \qquad  \left.+\left( 9\,{a}^{2}+3\,a-11\right) \,b+7\,{a}^{2}-11\,a+4\right)\\
T_4(a,b,c) &=& \left(\left( 4\,{d}^{2}+14\,d+7\right) \,{e}^{2}+\left( 14\,{d}^{2}+\left( \left( \left( -2\,a-8\right) \,b-8\,a-8\right) \,c+\left( -8\,a-8\right) \,b\right.\right.\right.\\
&& \qquad \left.\left.\left.-8\,a+32\right) \,d+\left( \left( -10\,a-16\right) \,b-16\,a-8\right) \,c+\left( -16\,a-8\right) \,b-8\,a+11\right) \,e\right.\\
&& \qquad \qquad  \left.+7\,{d}^{2}+\left( \left( \left( -10\,a-16\right) \,b-16\,a-8\right) \,c+\left( -16\,a-8\right) \,b-8\,a+11\right) \,d\right.\\
&& \qquad \qquad \qquad  \left.+\left( \left( 2\,a+4\right) \,{b}^{2}+\left( 2\,{a}^{2}+16\,a+10\right) \,b+4\,{a}^{2}+10\,a+3\right) \,{c}^{2}\right.\\
&& \qquad \qquad \qquad \qquad  \left.+\left( \left( 2\,{a}^{2}+16\,a+10\right) \,{b}^{2}+\left( 16\,{a}^{2}+16\,a-8\right) \,b+10\,{a}^{2}-8\,a\right.\right.\\
&& \qquad \qquad \qquad \qquad \qquad\left. \left.-5\right) \,c+\left( 4\,{a}^{2}+10\,a+3\right) \,{b}^{2}+\left( 10\,{a}^{2}-8\,a-5\right) \,b\right.\\ && \qquad \qquad \qquad \qquad \qquad \qquad \left.+3\,{a}^{2}-5\,a+2\right)
\end{eqnarray*}
and
\begin{eqnarray*}
T(a,b,c,d,e)&=&\left( 2\,{d}^{2}+2\,d\right) \,{e}^{2}+\left( 2\,{d}^{2}+\left( 2-8\,a\,b\,c\right) \,d-8\,a\,b\,c\right) \,e-8\,a\,b\,c\,d\\
&& \qquad +\left( \left( 3\,{a}^{2}+3\,a\right) \,{b}^{2}+\left( 3\,{a}^{2}+3\,a\right) \,b\right) \,{c}^{2}+\left( \left( 3\,{a}^{2}+3\,a\right) \,{b}^{2}+\left( 3\,{a}^{2}-5\,a\right) \,b\right) \,c \geq 0
 \end{eqnarray*}
 Then $z\, _3F_2(a,b,c;d,e;z)$ is in $KS^{*}$.
\ethm
\bpf The function $f(z)= z\, _3F_2(a,b,c;d,e;z)$ is defined by (\ref{eq1}), where $A_n$  is as in $(\ref{2f2eq1})$.\\

For $a,b,c > 0$, we observe that $de \geq 2\,a\,b\,c$ and  $$d+e \geq Max\left\{ a+b+c,\, \frac{1}{2}(ab+bc+ac+2(a+b+c)-1-2abc),2[ab+bc+ac]-3abc \right\}.$$ By Theorem \ref{thm12f2}, this condition implies that the sequence $\{nA_n\}$ is non-increasing. To prove $f$ is starlike. We need to show that the sequence $\{nA_n-(n+1)A_{n+1}\}$ is also non-increasing using Lemma \ref{lemeq5}. Let
 \begin{eqnarray*}
B_n = nA_n-(n+1)A_{n+1}\, \, \, {\ and}\, \, \, B_{n+1}&=& (n+1)A_{n+1} -  (n+2)A_{n+2}
 \end{eqnarray*}
Using $A_n$,  we find that
 \begin{eqnarray}\label{eqnP3}
B_{n}-B_{n+1}&=&nA_n -2(n+1)A_{n+1} +  (n+2)A_{n+2}\\
&=& A_n \left[n -2(n+1)\left(\frac{A_{n+1}}{A_n}\right) +  (n+2)\left(\frac{A_{n+2}}{A_n}\right)\right]\nonumber
 \end{eqnarray}
Where
 \begin{eqnarray*}
\frac{A_{n+1}}{A_n} = \frac{(a+n-1)(b+n-1)(c+n-1)}{(d+n-1)(e+n-1)n}
 \end{eqnarray*}
and
 \begin{eqnarray*}
\frac{A_{n+2}}{A_n} = \frac{(a+n)(a+n-1)(b+n)(b+n-1)(c+n)(c+n-1)}{(d+n)(d+n-1)(e+n)(e+n-1)n(n+1)}
 \end{eqnarray*}
After some simplification the equation (\ref{eqnP3}) implies
 \begin{eqnarray}\label{eqnP4}
B_{n}-B_{n+1}&=&\frac{ A_n P(n)}{n(n+1)(d+n)(d+n-1)(e+n)(e+n-1)}
 \end{eqnarray}
 \begin{eqnarray*}
P(n)&=&n^2(n+1)(d+n)(d+n-1)(e+n)(e+n-1)\nonumber\\
        &&\qquad  -2(n+1)^2(d+n)(e+n)(a+n-1)(b+n-1)(c+n-1)\nonumber\\
        && \qquad  \qquad +(n+2)(a+n)(a+n-1)(b+n)(b+n-1)(c+n)(c+n-1)\nonumber\\
    &=& { \left( e+d-c-b-a\right) \,\left( e+d-c-b-a+1\right)}\,n^5 \nonumber\\
    && +{2\,\left( e+d-c-b-a+1\right) \,\left( d\,e-b\,c-a\,c-a\,b+1\right)}\, n^4\nonumber\\
    && + \left(\left( {d}^{2}+d-1\right) \,{e}^{2}+\left( {d}^{2}+\left( \left( -2\,b-2\,a\right) \,c-2\,a\,b+1\right) \,d\right.\right.\nonumber\\
    &&\qquad \left.\left.+\left( \left( -2\,a-2\right) \,b-2\,a+4\right) \,c+\left( 4-2\,a\right) \,b+4\,a-3\right) \,e-{d}^{2}\right.\nonumber\\
    && \qquad \qquad \left.+\left( \left( \left( -2\,a-2\right) \,b-2\,a+4\right) \,c+\left( 4-2\,a\right) \,b+4\,a-3\right) \,d \right. \nonumber\\
    && \qquad \qquad  \qquad \left.+\left( {b}^{2}+\left( 4\,a+1\right) \,b+{a}^{2}+a-3\right) \,{c}^{2}+\left( \left( 4\,a+1\right) \,{b}^{2}+\left( 4\,{a}^{2}-9\right) \,b\right.\right.\nonumber\\
    && \qquad \qquad \qquad \qquad \left.\left.+{a}^{2}-9\,a+7\right) \,c+\left( {a}^{2}+a-3\right) \,{b}^{2}+\left( {a}^{2}-9\,a+7\right) \,b-3\,{a}^{2}\right.\nonumber\\
    && \qquad \qquad \qquad \qquad \qquad \left.+7\,a-4\right)\,{n}^{3}\nonumber\\
    && + \left(\left( {d}^{2}-d\right) \,{e}^{2}+\left( -{d}^{2}+\left( \left( \left( -2\,a-2\right) \,b-2\,a+4\right) \,c+\left( 4-2\,a\right) \,b+4\,a-3\right) \,d\right.\right.\nonumber \\
    && \qquad \left.\left.+\left( \left( 2-4\,a\right) \,b+2\,a\right) \,c+2\,a\,b-2\right) \,e+\left( \left( \left( 2-4\,a\right) \,b+2\,a\right) \,c+2\,a\,b-2\right) \,d\right.\nonumber \\
    && \qquad  \left.+\left( \left( 2\,a+1\right) \,{b}^{2}+\left( 2\,{a}^{2}+4\,a-5\right) \,b+{a}^{2}-5\,a+2\right) \,{c}^{2}+\left( \left( 2\,{a}^{2}+4\,a-5\right) \,{b}^{2}\right.\right.\nonumber \\
    && \qquad \qquad \left.\left.+\left( 4\,{a}^{2}-20\,a+11\right) \,b-5\,{a}^{2}+11\,a-4\right) \,c+\left( {a}^{2}-5\,a+2\right) \,{b}^{2}\right.\nonumber \\
    && \qquad \qquad \qquad \left.+\left( -5\,{a}^{2}+11\,a-4\right) \,b+2\,{a}^{2}-4\,a+2\right)\,{n}^{2}\nonumber\\
    &&+\left(\left( \left( \left( \left( 2-4\,a\right) \,b+2\,a\right) \,c+2\,a\,b-2\right) \,d+\left( \left( 2-2\,a\right) \,b+2\,a-2\right) \,c+\left( 2\,a-2\right) \,b\right.\right.\nonumber\\
     && \qquad \left.\left.-2\,a+2\right) \,e+\left( \left( \left( 2-2\,a\right) \,b+2\,a-2\right) \,c+\left( 2\,a-2\right) \,b-2\,a+2\right) \,d\right.\nonumber\\
     && \qquad \qquad \left.+\left( \left( {a}^{2}+3\,a-2\right) \,{b}^{2}+\left( 3\,{a}^{2}-7\,a+2\right) \,b-2\,{a}^{2}+2\,a\right) \,{c}^{2}\right.\nonumber\\
     && \qquad \qquad \qquad  \left.+\left( \left( 3\,{a}^{2}-7\,a+2\right) \,{b}^{2}+\left( -7\,{a}^{2}+11\,a-2\right) \,b+2\,{a}^{2}-2\,a\right) \,c\right.\nonumber\\
     && \qquad \qquad \qquad \qquad  \left.+\left( 2\,a-2\,{a}^{2}\right) \,{b}^{2}+\left( 2\,{a}^{2}-2\,a\right) \,b\right)\,n \nonumber\\
     && {-2\,\left( a-1\right) \,\left( b-1\right) \,\left( c-1\right) \,\left( d\,e-a\,b\,c\right)}
 \end{eqnarray*}
Our aim is to check that $P(n)$ is non-negative for all $n\geq1$. After few steps of calculation, we get
 \begin{eqnarray*}
T(n)&\geq & \left(\left( 5\,{d}^{2}+9\,d+2\right) \,{e}^{2}+\left( 9\,{d}^{2}+\left( \left( \left( -8\,a-8\right) \,b-8\,a\right) \,c-8\,a\,b+13\right) \,d\right.\right.\\
&& \qquad \left.\left.+\left( \left( -16\,a-8\right) \,b-8\,a\right) \,c-8\,a\,b+2\right) \,e+2\,{d}^{2}+\left( \left( \left( -16\,a-8\right) \,b-8\,a\right) \,c\right.\right.\\
&& \qquad \qquad \left.\left.-8\,a\,b+2\right) \,d+\left( \left( {a}^{2}+7\,a+3\right) \,{b}^{2}+\left( 7\,{a}^{2}+13\,a+3\right) \,b+3\,{a}^{2}+3\,a\right) \,{c}^{2}\right.\\
&& \qquad \qquad \qquad \left.+\left( \left( 7\,{a}^{2}+13\,a+3\right) \,{b}^{2}+\left( 13\,{a}^{2}-5\,a-5\right) \,b+3\,{a}^{2}-5\,a\right) \,c\right.\\
&& \qquad \qquad \qquad \qquad \left.+\left( 3\,{a}^{2}+3\,a\right) \,{b}^{2}+\left( 3\,{a}^{2}-5\,a\right) \,b\right)\,n\\
&& +\left( -3\,{d}^{2}-7\,d-2\right) \,{e}^{2}+\left( -7\,{d}^{2}+\left( \left( 8\,b+8\,a\right) \,c+8\,a\,b-11\right) \,d+\left( \left( 8\,a+8\right) \,b+8\,a\right) \,c\right.\\
&& \qquad \left.+8\,a\,b-2\right) \,e-2\,{d}^{2}+\left( \left( \left( 8\,a+8\right) \,b+8\,a\right) \,c+8\,a\,b-2\right) \,d+\left( \left( 2\,{a}^{2}-4\,a-3\right) \,{b}^{2}\right.\\
&& \qquad \qquad  \left.+\left( -4\,{a}^{2}-10\,a-3\right) \,b-3\,{a}^{2}-3\,a\right) \,{c}^{2}+\left( \left( -4\,{a}^{2}-10\,a-3\right) \,{b}^{2}\right.\\
&& \qquad \qquad \qquad \left.+\left( 5-10\,{a}^{2}\right) \,b-3\,{a}^{2}+5\,a\right) \,c+\left( -3\,{a}^{2}-3\,a\right) \,{b}^{2}+\left( 5\,a-3\,{a}^{2}\right) \,b
 \end{eqnarray*}
Put $n=1$, we have
 \begin{eqnarray*}
T(1)&\geq&\left( 2\,{d}^{2}+2\,d\right) \,{e}^{2}+\left( 2\,{d}^{2}+\left( 2-8\,a\,b\,c\right) \,d-8\,a\,b\,c\right) \,e-8\,a\,b\,c\,d\\
&& \qquad +\left( \left( 3\,{a}^{2}+3\,a\right) \,{b}^{2}+\left( 3\,{a}^{2}+3\,a\right) \,b\right) \,{c}^{2}+\left( \left( 3\,{a}^{2}+3\,a\right) \,{b}^{2}+\left( 3\,{a}^{2}-5\,a\right) \,b\right) \,c \geq 0
 \end{eqnarray*}
Hence by hypothesis,  the following inequalities holds true
$$P(n) \geq Q(n)\geq R(n)\geq S(n)\geq T(n)\geq T(1)\geq 0.$$
Therefore sequences $\{B_n\}$ and $\{nA_n-(n+1)A_{n+1}\}$,  is  non-increasing. We deduce that $f$ is starlike by Lemma \ref{lemeq5}. Also the function $f$ is close-to-convex with respect to $-\log(1-z)$. Since the conditions of Lemma \ref{lemeq2} are verified.
\epf

\bthm\label{thm42f2}
Let $a,b,\ and\ c > 0$. Suppose that  $$ d+e \geq max\{T_1(a,b,c),\,  T_2(a,b,c),\,  T_3(a,b,c)\}$$ where
 \begin{eqnarray*}
 T_1(a,b,c) &=&  \left( e+d-c-b-a+1\right) \,\left( e+d-c-b-a+2\right),\\
 T_2(a,b,c) &=& 2\,\left( e+d-c-b-a+2\right) \,\left( d\,e+2\,e+2\,d-b\,c-a\,c-c-a\,b-b-a+1\right),\\
  T_3(a,b,c)&=& \left(\left( {d}^{2}+7\,d+5\right) \,{e}^{2}+\left( 7\,{d}^{2}+\left( \left( -2\,b-2\,a-4\right) \,c+\left( -2\,a-4\right) \,b-4\,a+21\right) \,d\right.\right.\\
    && \qquad \left.\left.+\left( \left( -2\,a-6\right) \,b-6\,a-4\right) \,c+\left( -6\,a-4\right) \,b-4\,a+9\right) \,e+5\,{d}^{2}\right.\\
        && \qquad \qquad\left.+\left( \left( \left( -2\,a-6\right) \,b-6\,a-4\right) \,c+\left( -6\,a-4\right) \,b-4\,a+9\right) \,d\right.\\
                &&\qquad \qquad \qquad  \left.+\left( {b}^{2}+\left( 4\,a+3\right) \,b+{a}^{2}+3\,a+1\right) \,{c}^{2}+\left( \left( 4\,a+3\right) \,{b}^{2}\right.\right.\\
                    &&\qquad \qquad \qquad  \qquad \left.\left.+\left( 4\,{a}^{2}+4\,a-5\right) \,b+3\,{a}^{2}-5\,a-3\right) \,c+\left( {a}^{2}+3\,a+1\right) \,{b}^{2}\right.\\
                        &&\qquad \qquad \qquad \qquad \qquad  \left.+\left( 3\,{a}^{2}-5\,a-3\right) \,b+{a}^{2}-3\,a+2\right)
  \end{eqnarray*}
and satisfies
\begin{eqnarray*}
T(a,b,c,d,e)&=&\left( 2\,{d}^{2}+2\,d\right) \,{e}^{2}+\left( 2\,{d}^{2}+\left( 2-4\,a\,b\,c\right) \,d-4\,a\,b\,c\right) \,e-4\,a\,b\,c\,d\\
    && \qquad +\left( \left( {a}^{2}+a\right) \,{b}^{2}+\left( {a}^{2}+a\right) \,b\right) \,{c}^{2}+\left( \left( {a}^{2}+a\right) \,{b}^{2}+\left( {a}^{2}-3\,a\right) \,b\right) \,c\geq 0
 \end{eqnarray*}
then the Alexander transform is defined by (\ref{alexeq1}) is in $KS^{*}$.
\ethm
\bpf
Let $f(z) = z\, _3F_2(a,b,c;d,e;z)$, then from the definition of hypergeometric function $_3F_2$, we have
 \begin{eqnarray*}
f(z)=\sum_{n=1}^{\infty} \frac{(a)_{n-1}(b)_{n-1}(c)_{n-1}}{(d)_{n-1}(e)_{n-1}(1)_{n-1}}z^n
 \end{eqnarray*}
so that the corresponding Alexander transform defined by (\ref{alexeq1}) takes the form
 \begin{eqnarray*}
\Lambda_f(z)=\sum_{n=1}^{\infty}A_nz^n,
 \end{eqnarray*}
with\ $A_1 = 1$ and
 \begin{eqnarray}
A_n = \frac{(a)_{n-1}(b)_{n-1}(c)_{n-1}}{n(d)_{n-1}(e)_{n-1}(1)_{n-1}}, for\ n\geq 2.
 \end{eqnarray}
Using the definition of the shifted factorial notation, we have
 \begin{eqnarray*}
A_{n+1} &=& \frac{(a)_{n}(b)_{n}(c)_{n}}{(n+1)(d)_{n}(e)_{n}(1)_{n}}
 \end{eqnarray*}
and
 \begin{eqnarray*}
(n+1)A_{n+1}&=& \frac{(a+n-1)(b+n-1)(c+n-1)}{(d+n-1)(e+n-1)} A_n\\
 \end{eqnarray*}
After some simplification,
 \begin{eqnarray*}
nA_n - (n+1)A_{n+1} &=& nA_n -  \frac{(a+n-1)(b+n-1)(c+n-1)}{(d+n-1)(e+n-1)} A_n\\
&=& A_n \left[\frac{n(d+n-1)(e+n-1)-(a+n-1)(b+n-1)(c+n-1)}{(d+n-1)(e+n-1)} \right]\\
&=& \frac{A_n U(n) }{(d+n-1)(e+n-1)}
 \end{eqnarray*}
where
 \begin{eqnarray*}
 U(n) &=& n(d+n-1)(e+n-1)-(a+n-1)(b+n-1)(c+n-1)\\
 &\geq& {\left( e+d-c-b-a+1\right)} \,{n}^{2}\\
 && \qquad +\left( \left( d-1\right) \,e-d+\left( -b-a+2\right) \,c+\left( 2-a\right) \,b+2\,a-2\right) \,n\\
 && \qquad \qquad +\left( \left( 1-a\right) \,b+a-1\right) \,c+\left( a-1\right) \,b-a+1\\
 &=&{\left(d\,e+e+d-b\,c-a\,c-a\,b\right)} \,n-e-d+\left( \left( 1-a\right) \,b+a\right) \,c+a\,b
 \end{eqnarray*}
By hypothesis, $d+e\geq a\,b+b\,c+a\,c-abc$ and  so for every  $n\geq 1$, $$ U(n)\geq U(1)=  de-abc \geq 0$$ and since $de\geq abc$, we proved that the above condition satisfied. Thus the sequence $\{nA_n\}$ is decreasing. Next, we prove that $\{nA_n-(n+1)A_{n+1}\}$ is also decreasing.\\

Let $B_n = nA_n-(n+1)A_{n+1}$. Then after some manipulation, we get
 \begin{eqnarray*}
B_n-B_{n+1}&=&\frac{A_n U(n)}{(d+n-1)(e+n-1)}-\frac{A_{n+1} U(n+1)}{(d+n)(e+n)}\\
&=&\frac{A_n U(n)}{(d+n-1)(e+n-1)}-\frac{A_{n} U(n+1)}{(d+n)(e+n)(n+1)}\\ && \qquad \qquad \times\left( \frac{(a+n-1)(b+n-1)(c+n-1)}{(d+n-1)(e+n-1)}\right)\\
&=&\frac{A_n C(n)}{(d+n-1)(c+n-1)(d+n)(e+n)(n+1)}
 \end{eqnarray*}
where
 \begin{eqnarray*}
C(n)&=&U(n)(d+n)(e+n)(n+1)-U(n+1)(a+n-1)(b+n-1)(c+n-1)\\
&=& {\left( e+d-c-b-a+1\right) \,\left( e+d-c-b-a+2\right)} \,{n}^{4}\\
&& +{\left( e+d-c-b-a+2\right) \,\left( d\,e-b\,c-a\,c+c-a\,b+b+a-1\right)}\,{n}^{3}\\
&& +{2\,\left( e+d-c-b-a+2\right) \,\left( d\,e-b\,c-a\,c+c-a\,b+b+a-1\right)}\,{n}^{2}\\
&& + \left(\left( {d}^{2}-d\right) \,{e}^{2}+\left( -{d}^{2}+\left( \left( 2-2\,a\,b\right) \,c+2\,b+2\,a-3\right) \,d\right.\right.\\
    &&\qquad  \left.\left.+\left( \left( 2-2\,a\right) \,b+2\,a-2\right) \,c+\left( 2\,a-2\right) \,b-2\,a+2\right) \,e\right.
 \end{eqnarray*}
 \begin{eqnarray*}
        && \qquad \qquad \left. +\left( \left( \left( 2-2\,a\right) \,b+2\,a-2\right) \,c+\left( 2\,a-2\right) \,b-2\,a+2\right) \,d\right.\\
            && \qquad \qquad \qquad\left.+\left( \left( 2\,a-1\right) \,{b}^{2}+\left( 2\,{a}^{2}-4\,a+1\right) \,b-{a}^{2}+a\right) \,{c}^{2}\right.\\
                && \qquad \qquad \qquad \qquad \left.+\left( \left( 2\,{a}^{2}-4\,a+1\right) \,{b}^{2}+\left( -4\,{a}^{2}+6\,a-1\right) \,b+{a}^{2}-a\right) \,c\right.\\
                     && \qquad \qquad \qquad  \qquad \qquad \left.+\left( a-{a}^{2}\right) \,{b}^{2}+\left( {a}^{2}-a\right) \,b\right)\,n\\
                        && \qquad \qquad \qquad  \qquad \qquad \qquad -{\left( a-1\right) \,\left( b-1\right) \,\left( c-1\right) \,\left( 2\,d\,e-a\,b\,c\right)}
 \end{eqnarray*}
Our aim is to check that $C(n)$ is non-negative for all $n\geq1$. After manipulation,  we get
 \begin{eqnarray*}
F(n)&\geq&\left(\left( 3\,{d}^{2}+7\,d+2\right) \,{e}^{2}+\left( 7\,{d}^{2}+\left( \left( \left( -2\,a-4\right) \,b-4\,a\right) \,c-4\,a\,b+11\right) \,d\right.\right.\\
    && \qquad \left. \left.+\left( \left( -6\,a-4\right) \,b-4\,a\right) \,c-4\,a\,b+2\right) \,e+2\,{d}^{2}+\left( \left( \left( -6\,a-4\right) \,b-4\,a\right) \,c \right.\right.\\
        && \qquad \qquad \left.\left.-4\,a\,b+2\right) \,d+\left( \left( 2\,a+1\right) \,{b}^{2}+\left( 2\,{a}^{2}+4\,a+1\right) \,b+{a}^{2}+a\right) \,{c}^{2}\right.\\
            && \qquad \qquad \qquad \left.+\left( \left( 2\,{a}^{2}+4\,a+1\right) \,{b}^{2}+\left( 4\,{a}^{2}-4\,a-3\right) \,b+{a}^{2}-3\,a\right) \,c\right.\\
                && \qquad \qquad \qquad \qquad \left.+\left( {a}^{2}+a\right) \,{b}^{2}+\left( {a}^{2}-3\,a\right) \,b\right)\,n\\
&& +\left( -{d}^{2}-5\,d-2\right) \,{e}^{2}+\left( -5\,{d}^{2}+\left( \left( \left( 4-2\,a\right) \,b+4\,a\right) \,c+4\,a\,b-9\right) \,d\right.\\
    && \qquad  \left.+\left( \left( 2\,a+4\right) \,b+4\,a\right) \,c+4\,a\,b-2\right) \,e-2\,{d}^{2}+\left( \left( \left( 2\,a+4\right) \,b+4\,a\right) \,c\right.\\
        && \qquad \qquad \left.+4\,a\,b-2\right) \,d+\left( \left( {a}^{2}-a-1\right) \,{b}^{2}+\left( -{a}^{2}-3\,a-1\right) \,b-{a}^{2}-a\right) \,{c}^{2}\\
            && \qquad \qquad \qquad +\left( \left( -{a}^{2}-3\,a-1\right) \,{b}^{2}+\left( -3\,{a}^{2}+a+3\right) \,b-{a}^{2}+3\,a\right) \,c\\
                && \qquad \qquad \qquad \qquad +\left( -{a}^{2}-a\right) \,{b}^{2}+\left( 3\,a-{a}^{2}\right) \,b
 \end{eqnarray*}
Put $n=1$ in the above, we get
 \begin{eqnarray*}
F(1)&=&\left( 2\,{d}^{2}+2\,d\right) \,{e}^{2}+\left( 2\,{d}^{2}+\left( 2-4\,a\,b\,c\right) \,d-4\,a\,b\,c\right) \,e-4\,a\,b\,c\,d\\
    && \qquad +\left( \left( {a}^{2}+a\right) \,{b}^{2}+\left( {a}^{2}+a\right) \,b\right) \,{c}^{2}+\left( \left( {a}^{2}+a\right) \,{b}^{2}+\left( {a}^{2}-3\,a\right) \,b\right) \,c\geq 0
 \end{eqnarray*}

For every $n \geq 1$, $C(n)\geq D(n)\geq E(n) \geq F(n)\geq F(1)$. Which is true \\

Therefore, the sequences $\{B_n\}$ and  $\{nA_n - (n+1)A_{n+1}\}$ are non-increasing. The Alexander transform $\Lambda_f(z)$ is starlike in unit disc $\D$ by Lemma \ref{lemeq5} and already verified  $\{nA_n\}$ is non-increasing. Therefore, the Alexander transform $\, \Lambda_f(z)$  is close-to-convex with respect to $-\log(1-z)$ using Lemma \ref{lemeq2}. Thus we have $\Lambda_f(z) \in KS^{*}$. Which completes the proof of Theorem.
\epf

\end{document}